\documentclass[11pt]{article}
\title {A volume-ish theorem for the Jones polynomial\\of alternating knots}
\author {Oliver T. Dasbach
\thanks{e-mail: kasten@math.lsu.edu, 
http://www.math.lsu.edu/$\sim$\!\! kasten}
\thanks{supported in part by an NSF grant }
\\Louisiana State University\\
Department of Mathematics\\Baton Rouge, LA 70803
\and
Xiao-Song Lin
\thanks{e-mail: xl@math.ucr.edu, http://www.math.ucr.edu/$\sim$\!\! xl}
\thanks{supported in part by an NSF grant }
 \\ University of California,
    Riverside \\Department of Mathematics\\Riverside, CA 92521 - 0135
}

\date{}

\usepackage{amsmath,amssymb,amsthm, graphics}

\newtheorem{theorem}{Theorem}[section]

\newtheorem{proposition}[theorem]{Proposition}

\newtheorem{corollary}[theorem]{Corollary}

\newcommand {\order}[1]{\vert #1 \vert}

\begin{document} 
\maketitle
\begin{abstract} 
The Volume conjecture claims that the hyperbolic Volume of a knot is 
determined by the colored Jones polynomial.

The purpose of this article is to show a Volume-ish theorem for alternating 
knots in terms of the Jones polynomial, rather than the colored Jones 
polynomial: The ratio of the Volume and certain sums of coefficients of the 
Jones polynomial is bounded from above and from below by constants.

Furthermore, we give experimental data on the relation of the growths of the 
hyperbolic volume and the coefficients of the Jones polynomial, both
for alternating and non-alternating knots. 
\end {abstract}

\section {Introduction}

Since the introduction of the Jones polynomial, there is a strong desire
to have a geometrical or topological interpretation for it rather than a
combinatorial definition.

The, arguably, first major success in this direction was the proof of the 
Melvin-Morton Conjecture by Rozanski, Bar-Natan and Garoufalidis (see 
\cite{BN-G:Melvin-Morton} and compare with 
\cite{Vaintrob, Chmutov1, LW:RandomWalk} 
for different proofs): 
The Alexander polynomial is determined by the so-called 
colored Jones polynomial. For a knot $K$ the colored Jones polynomial  
is given by the Jones polynomial and the Jones polynomials of cablings
of $K$. 

The next major conjecture that relates the Jones polynomial and its 
offsprings to classical topology and geometry was the Volume conjecture of 
Kashaev, Murakami and Murakami (e.g. \cite{Murakami:VolumeConjecture}).
This conjecture states that the colored Jones polynomial determines the
Gromov norm of the knot complement. For hyperbolic knots the Gromov norm is
proportional to the hyperbolic volume.  

A proof of the Volume conjecture for all knots would also imply that the 
colored Jones polynomial detects the unknot 
\cite{Murakami:VolumeConjecture}.  This problem is still wide open; 
even for the Jones polynomial there is no counterexample known 
(see e.g. \cite{DH:DoesThe}).

The purpose of this paper is to show a relation
of the coefficients of the Jones polynomial and the hyperbolic volume of 
alternating knot complements. More specifically we prove:

\medskip

\noindent 
{\bf Volume-ish Theorem. } \sl
For an alternating, prime, non-torus knot 
$K$ let $$V_K(t)= a_n t^n + \dots + a_m t^m$$ be the Jones 
polynomial of $K$. Then 

$$ 2 v_0 (\max(\vert a_{m-1} \vert, \vert a_{n+1} \vert) -1)  \leq 
\mbox{Vol}(S^3-K) \leq 10 v_0 (|a_{n+1}|+|a_{m-1}|-1). $$
Here, $v_0 \approx 1.0149416$ is the Volume of an ideal regular hyperbolic 
tetrahedron.
\em
\medskip

For the proof of this theorem we make use of Marc Lackenby's result 
and its proof that the hyperbolic volume is bounded from above and below 
linearily by the twist number. 
His bound was improved by Ian Agol and Dylan Thurston.

In an appendix we give some numerical data on the relation of other 
coefficients and the hyperbolic volume, both for alternating and for 
non-alternating knots. These data gives some hope for a 
Volume-ish theorem for non-alternating knots as well.

\bigskip 
\noindent
{\bf Acknowledgement: } The first author would like to thank Ian Agol, 
Joan Birman,  Charlie Frohman, Vaughan Jones, Lou Kauffman, James Oxley, 
and Neal Stoltzfus for helpful conversations at various occasions. 

\section {The Jones Polynomial evaluation of the Tutte Polynomial}

Our goal is to relate the hyperbolic volume of alternating knot 
complements to the coefficients of the Jones polynomial.
For this we will make use of the computation of the Jones polynomial of
alternating links via the Tutte polynomial.

\subsection{Notations}

First, we need a few notations. Our objects are multi-graphs, i.e. parallel 
edges are allowed. Two edges are called parallel if they connect the same 
two vertices. 

\medskip

\noindent
{\bf Notation} 
\begin{enumerate}
\item A multi graph $G=(V,E)$ has a set $V$ of vertices and a set $E$ of edges.
\item We denote by $\tilde G$ a spanning subgraph of $G$, 
where parallel edges are deleted. See Figure \ref{spanning_subgraph}. 
The subgraph $\tilde G=(V, \tilde E)$ has a set of vertices $V$ and a set $\tilde E$ of edges.
\item Since $\tilde G$ is a subgraph of $G$, each edge $e \in \tilde E$ in 
$\tilde G$ comes with a multiplicity $\mu(e)$. 
This is the number of edges in $G$ 
that are parallel to $e$. For example, the graph in Figure 
\ref{spanning_subgraph} has one edge with multiplicity $2$, one 
edge with multiplicity $3$, one edge with multiplicity $4$ and all other edges have multiplicity $1$. 
It is convenient for us to define $n(j)$ to be the
number of edges $e \in \tilde E$ with $\mu(j) \geq j$.
In particular $n(1)=|\tilde E|$.
The graph in Figure \ref{spanning_subgraph} has: $n(2)=3, n(3)=2, n(4)=1.$

\begin{figure}[h]
\centerline{{\scalebox{0.55}
{\includegraphics{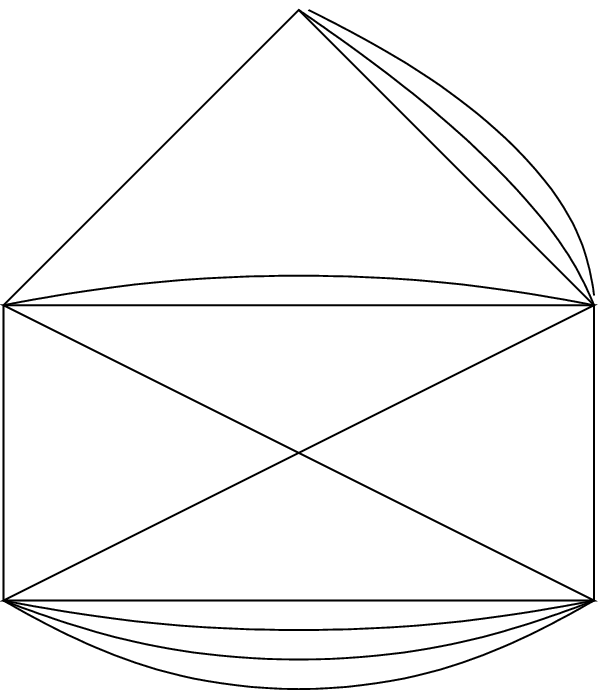}} \qquad
\scalebox{0.55} {\includegraphics{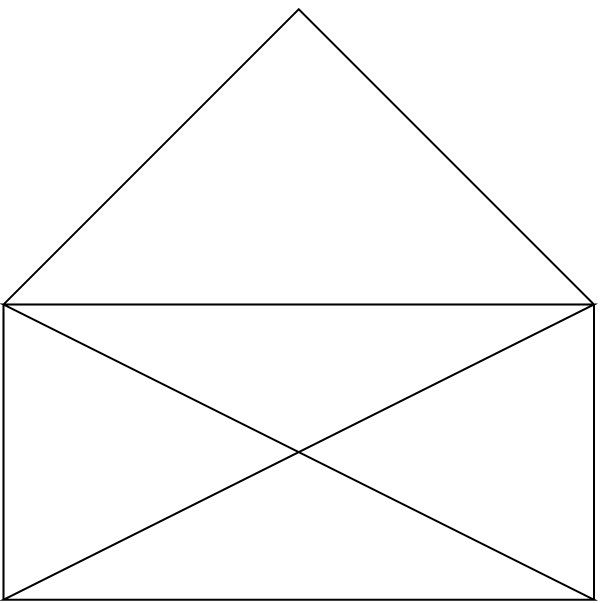}}}
}
\caption{A multi-graph $G$ and its spanning subgraph $\tilde G$} 
\label{spanning_subgraph}
\end{figure}

\item The number of components of a graph $G$ is $k(G)$. If $V$ is apparend 
from the context and
$G=(V,E)$ we set $k(E):=k(G)$.
\item The Tutte polynomial of a multi-graph $G$ (see e.g. \cite{Bollobas:Book})
\begin{eqnarray*}
T_G(x,y)&:=&\sum_{F\subseteq E} 
(x-1)^{k(F)-k(E)}(y-1)^{\vert F \vert - \vert V \vert + k(F)}
\end{eqnarray*}
\end{enumerate}

\subsubsection{The Tutte polynomial and the Jones polynomial for alternating
links}

Let $K$ be an alternating link with an alternating plane projection $P(K)$.
The region of the projection can be colored with 2 colors, say, purple and 
gold, such that two adjacent faces have different colors.

Two graphs are assigned to the projection, one corresponding to
the purple colored regions and one to the golden regions. 
Every region gives rise to a vertex in the graph and two vertices are
connected by an edge if the corresponding regions are adjacent to a common
crossing. These graphs are called checkerboard graphs.

\begin{figure}[h]
\centerline{{\scalebox{1.15}
{\includegraphics{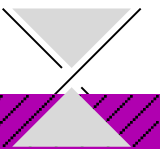}} \qquad
\scalebox{1.15} {\includegraphics{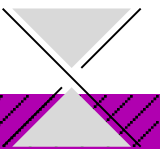}}}
}
\caption{A positive or negative sign for the shaded region in the checkerboard graph} 
\label{Figure_sign}
\end{figure}

Each edge comes with a sign as in Figure \ref{Figure_sign}. 

For an alternating link $K$ all edges are either positive or negative.
Thus we have a positive checkerboard graph and a negative checkerboard graph.
These two graphs are dual to each other.
Let $G$ be the positive checkerboard graph, $a$ be the number of vertices
in $G$ and $b$ be the number of vertices in the negative checkerboard graph.

It is well known, e.g.  \cite{Bollobas:Book}, 
that the Jones Polynomial of an alternating link $K$ with positive checkerboard
 graph $G$ satisfies:
$$V_K(t)=(-1)^w t^{(b-a+3 w)/4} T_G(-t,-1/t).$$

Here, $w$ is the writhe number, i.e. the algebraic crossing number of the link 
projection.

Since we are interested in the absolute values of coefficients of the Jones polynomial, all information relevant to us is contained in the evaluation 
$T_G(-t,-1/t)$ of the Tutte polynomial.

\subsection{Reduction of multiple edges to simple edges}

Our first step is to reduce the computation of the Tutte polynomial of a 
multi-graph to the computation of a weighted Tutte polynomial of a spanning 
simple graph.

Again, for a graph $G=(V,E)$ we denote by $\tilde G=(V,\tilde E)$ a spanning
simple graph. If $G$ is connected graph $G$ without vertices of valency $1$, 
i.e. without loops then:

\begin{eqnarray*}
T_G(-t,-1/t)&=&\sum_{F \subseteq E}(-t-1)^{k(F)-1} \left (-\frac 1 t-1
  \right)^{\vert F
  \vert-\vert V \vert+k(F)}\\
&=&\sum_{\tilde F\subseteq \tilde E} (-t-1)^{k(\tilde F)-1}
  \left(-\frac 1 t -1 \right )^{-\vert V \vert + k(\tilde F)} \times\\
&&
\left ( 
\sum_{\stackrel {r(e_1)=1, \dots, r(e_j)=1}{e_1,\dots,e_j \in \tilde
  F}}^{\mu(e_1),\dots,\mu(e_j)} {\mu(e_1) \choose r(e_1)} \cdots
  {\mu(e_j) \choose r(e_j)} \left (- \frac 1 t -1
  \right)^{r(e_1)+\dots+r(e_j)} \right )\\
&=& \sum_{\tilde F\subseteq \tilde E} \left ( (-t-1)^{k(\tilde F)-1}
  \left(-\frac 1 t -1 \right )^{-\vert V \vert + k(\tilde F)}
\prod _{e \in \tilde F} \left ( \left (-\frac 1 t \right)^{\mu(e)}
- 1 \right ) \right ).
\end{eqnarray*}

So, with
\begin{eqnarray*}
P(m)&:=& \frac {\left ( \left ( -\frac 1 t \right)^{m} - 1 \right )} 
{\left (-\frac 1 t -1 \right )}\\
&=& 1 - t^{-1} + t^{-2} - \dots \pm t^{-m+1}
\end{eqnarray*}
we have

\begin{eqnarray} \label{weighted_Tutte}
T_G(-t,-1/t)&=&\sum_{\tilde F\subseteq \tilde E} \left ( (-t-1)^{k(\tilde F)-1}
  \left(-\frac 1 t -1 \right )^{\vert \tilde F \vert 
-\vert V \vert + k(\tilde F)}
\prod _{e \in \tilde F} P(\mu(e))   \right ).
\end{eqnarray}

\subsection {The highest terms of $T_G(-t,-1/t)$}

\begin{proposition} \label{coefficients Tutte}
Let $G=(V,E)$ be a planar multi-graph with spanning simple graph
$\tilde G=(V,\tilde E)$. Let the Tutte Polynomial evaluation 
$$T_G(-t,-1/t)= a_n t^n + a_{n+1} t^{n+1}+ \dots + a_{m-1} t^{m-1} + a_m t^m,$$
for suitable $n$ and $m$.

Then the  coefficients of the highest degree terms of $T_G(-t,-1/t)$ are:

\begin{enumerate}
\item The highest degree term $t^m$ of $T_G(-t,-1/t)$ in $t$ is 
$t^{\order V -1}$ with coefficient $$a_m=(-1)^{\order V -1}.$$

\item The second highest degree term is $t^{\order V -2}$ with coefficient:
$$a_{m-1}=(-1)^{\order V-1} \left (\order V -1 - \order {\tilde E}\right ).$$

Note that $\vert a_{m-1}\vert = \order {\tilde E}+1 - \order V$.

\item The third highest degree term is $t^{\order V -3}$ with coefficient:

$$(-1)^{\order V}\left ( - {\order V -1 \choose 2} + (\order V -2)
  \order{\tilde E}- n(2) - {\order{\tilde E} \choose 2}+ \mbox{tri}
\right ),$$

where $\mbox{tri}$ is the number of triangles in $\tilde E$.

This term equals
$$a_{m-2}=(-1)^{\order V} \left ( - {|a_{m-1}|+1 \choose 2} -n(2)
  +\mbox{tri} \right )$$

\end{enumerate}
\end{proposition}

\begin{proof}
It is easy to see that 
$\vert \tilde F \vert - \vert V \vert +k(F) \geq 0$ for all $F$.

Therefore, 
$$
\left (-\frac 1 t -1 \right )^{\vert \tilde F \vert -  \vert V \vert +k(F)} 
\prod _{e \in \tilde F} P(\mu(e)) =
\pm 1 + \mbox{ higher terms in } t^{-1}$$

This means to determine the highest terms of $T_G(-t,-1/t)$  
we have to analyze terms where $k(\tilde F)$ is large.

\begin{enumerate}
\item[-] {\bf Case: } $k(\tilde F)=\order{V}$.

This means: $\order{\tilde F}=0.$
Thus the contribution in the sum in (\ref{weighted_Tutte}) is: 

$$(-t-1)^{\order V-1} = (-1)^{\order V -1} \left ( t^{\order V-1}+
(\order V -1) \, t^{\order V-2} + {\order V -1 \choose 2} \, t^{\order
  V-3} + \dots +1 \right )$$ 

\item[-] {\bf Case: } $k(\tilde F)=\order{V}-1$.

This means: $\order {\tilde F}=1.$
Thus the contribution is:

$$ (-t-1)^{\order V -2} \, \sum_{e \in \tilde E} P(\mu(e))$$

Recall our notation:  For $G$ a multi-graph $n(j)$ is the number 
of edges in $\tilde E$ of multiplicity $\geq j$.
In particular $n(1)=\order {\tilde E}$ 

Thus
$$ \sum_{e \in \tilde E} P(\mu(e))= \order {\tilde E} - n(2) t^{-1} +
n(3) t^{-2} - n(4) t^{-3} + \dots.$$

\item[-] {\bf Case: } $k(\tilde F)=\order{V}-2$.

In this case either $\order{\tilde F}=2$ or $\order {\tilde F}=3$ and
$\tilde F$ is a triangle. Thus the contribution is:
\begin{eqnarray*}
&&\sum_{e, f \in \tilde E} (-t-1)^{\order V -3} \, P(\mu(e)) \,
P(\mu(f))+\\ 
&&\sum_{\stackrel{e,f,g \in \tilde E}{ (e,f,g) \, \mbox{\tiny
      triangle}}} (-t-1)^{\order V-3} \left (-\frac 1 t-1 \right ) 
P(\mu(e)) P(\mu(f)) P(\mu(g))
\end{eqnarray*}
\end{enumerate}

By combining these computations we get the result.
\end{proof}

\section {An Algebraic Point of View}

It is interesting to formulate the results of 
Proposition \ref{coefficients Tutte} in a purely algebraic way:

Let $G$ be a multi-graph and $A$ it's $N \times N$ adjacency matrix.
Define a matrix $\tilde A$ to be the matrix that we get from $A$ by replacing
every non-zero entry $A$ by $1$. Thus, $\tilde A$ has only $1$ and $0$ as entries. Furthermore, $n(2)$ is half the number of entries in $A$ that are $\geq 2$.

Then:

\begin{corollary} Let
$$T_G(-t,-1/t)= a_n t^n + a_{n+1} t^{n+1}+ \dots + a_{m-1} t^{m-1} + a_m t^m,$$
be the Jones evaluation of the Tutte Polynomial of a planar graph $G$.

Then:
\begin {eqnarray*}
|a_m|&=&1\\
|a_{m-1}|&=& \frac 1 2 \mbox{ trace } {\tilde A}^2 -1 - N\\
|a_{m-2}|&=& { |a_{m-1}|+1 \choose 2} + n(2) - \frac 1 6 \mbox { trace } \tilde A^3 
\end {eqnarray*}
\end{corollary}

\begin{proof}
This is an immediate consequence of Proposition \ref{coefficients Tutte}
together with the well-known results that for an adjacency matrix $A$
the number of edges is $\frac 1 2 \mbox{ trace }A^2$ and the number of 
triangles is $\frac 1 6 \mbox{ trace } A^3$ (see e.g \cite{Biggs:AlgebraicGraphTheory}). 
\end{proof}

\section {Twist number and the Volume of an hyperbolic alternating
  knot}

For a reference on hyperbolic structures on knot complements we refer e.g.
to \cite{CR:HyperbolicStructures}. 
The figure eight knot has minimal volume among all hyperbolic 
knot complements \cite{CM:MinimalVolume}.
Moreover, for a hyperbolic knot $K$ with crossing number $c>4$ 
by a result of Colin Adams (see \cite{CR:HyperbolicStructures}), 
the hyperbolic volume of its complement is bounded from above by  
$$Vol(S^3-K) \leq (4 c -16) v_0,$$ where $v_0$ is the volume of a
regular ideal hyperbolic tetrahedron.

For alternating knot complements a better general upper bound is known in
terms of the twist number. As shown by Bill Menasco 
\cite{Menasco:AlternatingHyperbolic} a non-torus alternating knot is 
hyperbolic.

\begin{figure}
\centerline{{\scalebox{0.38}
{\includegraphics{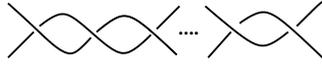}}} \label{twist}
}
\caption{A twist in a diagram of a knot} 
\end{figure}

The twist number of a diagram of an alternating knot 
is the minimal number of twists (see Figure 
\ref{twist}) in it. Here, a twist can consist of a single crossing.
For example, the diagram of the knot in Figure \ref{SampleKnot} has twist 
number 8. 

A twist corresponds to parallel edges in one of the checkerboard graphs.
Let $D$ be a diagram for an alternating knot $K$, $G=(V,E)$ be one of the two 
checkerboard graphs and let $G^*=(V^*,E^*)$ be the other, which is dual to 
$G$. In particular, $|E|=|E^*|$. As above we denote by $\tilde G=(V, \tilde E)$
the reduced (simple) graph and by $\tilde G^*=(V^*, \tilde E^*)$ the reduced
graph of the dual graph. 

We can now define the twist number by

\begin{eqnarray*}
T(K)&:=&|E|-(|E|-|\tilde E|)-(|E^*|-|\tilde E^*|)\\
&=& |E|-(|E|-|\tilde E|)-(|E|-|\tilde E^*|)\\
&=& |\tilde E|+|\tilde E^*|-|E|.
\end{eqnarray*}

It is an easy exercise to see that 
\begin{enumerate}
\item $T(K)$ is indeed realized as the twist number of a diagram of $K$.
\item $T(K)$ is an invariant of all alternating projections of $K$.
This follows from the Tait-Menasco-Thistlethwaite flyping theorem
\cite{MT:FlypingConjecture}. Below we will give a different argument for it.
\end{enumerate}

\begin {theorem}[Lackenby \cite{Lackenby:Volume}, Agol, D. Thurston]
\label{Lackenby}
$$v_0  \, (T(K)-2) \leq Vol(S^3-K) <10 v_0 (T(K)-1),$$
where $Vol(S^3-K)$ is the hyperbolic Volume and 
$v_0$ is the volume of an ideal regular hyperbolic tetrahedron.
\end{theorem}

\section {Coefficients of the Jones polynomial}

For an alternating knot $K$ with reduced alternating diagram $D$ with
$c$ crossings it holds \cite{Thistlethwaite:SpanningTreeExpansion, Kauffman:StateModels, Murasugi:Alternating}:

\begin{enumerate}
\item The span of the Jones polynomial is $c$.
\item The signs of the coefficients are alternating
\item The absolute values of the highest and lowest coefficients are $1$.
\end{enumerate}

Proposition \ref{coefficients Tutte} immediately leads to:

\begin{theorem} \label{TwistnumberinJones}
Let $V_K(t)=a_n t^n + a_{n+1} t^{n+1} + \dots + a_m t^m$ be the Jones
polynomial of an alternating knot $K$ and let $G=(V,E)$ be a checkerboard
graph of a reduced alternating projection of $K$.

Then

\begin{enumerate}
\item $\order {a_n}= \order{a_m}=1$
\item $\order {a_{n+1}}+\order{a_{m-1}}=T(K)$
\item $$\order {a_{n+2}}+\order{a_{m-2}}+ \order{a_{m-1}}\order{a_{n+1}}=
\frac {T(K)+T(K)^2} 2 +n(2) + n^*(2) -\mbox {tri}-\mbox{tri}^*,$$

where $n(2)$ is the number of edges in $\tilde E$ of multiplicity $>1$ and 
$n^*(2)$ the corresponding number in the dual checkerboard graph.

The number $\mbox{tri}$ is the number of triangles in the graph 
$\tilde G=(V,\tilde E)$
and $\mbox{tri}^*$ corresponds to $\mbox{tri}$ in the dual graph.
\item In particular, the twist number is an invariant of reduced alternating 
projections of the knot.
\end{enumerate}

\end{theorem}

\begin{proof}
Let $K$ be an alternating knot with a reduced alternating
 diagram $D$ having a checkerboard graph $G=(V,E)$ and dual graph
 $G^*=(V^*, E^*)$.

We have: 
$$\order E = \order {E^*} \quad \mbox{ and } \quad 
\order V+\order {V^*}=\order E+2.$$

Recall, that the twist number $T(K)$ of $K$ is defined as:

\begin{eqnarray*} 
T(K)&:=& \order E - (\order E - \order {\tilde E}) - (\order E-\order
{\tilde {E^*}})\\
&=& \order {\tilde E} + \order {\tilde {E^*}}-\order E\\
&=&  (\order {\tilde E} - \order V+1) + (\order {\tilde {E^*}}-\order {V^*}+1)
\end{eqnarray*}

The identities now follow from Proposition \ref{coefficients Tutte}.
\end{proof}

\noindent 
{\bf Volume-ish Theorem.} \sl
For an alternating, prime, non-torus knot 
$K$ let $$V_K(t)= a_n t^n + \dots + a_m t^m$$ be the Jones 
polynomial of $K$. Then 

$$ 2 v_0 (\max(\vert a_{m-1} \vert, \vert a_{n+1} \vert) -1)  \leq 
\mbox{Vol}(S^3-K) \leq 10 v_0 (|a_{n+1}|+|a_{m-1}|-1). $$

Here, $v_0 \approx 1.0149416$ is the Volume of an ideal regular hyperbolic 
tetrahedron.
\em

\medskip

\begin{proof}
The upper bound follows from Theorem \ref{Lackenby} and 
Theorem \ref{TwistnumberinJones}. 
For the lower bound we need a closer look at \cite{Lackenby:Volume}.

We can suppose that $K$ admits a diagram such that both checkerboard graphs
are inbedded so that every pair of edges connecting the same two 
vertices are adjacent to each other in the plane. 
This can be done by flypes.

Suppose $G_p=(V_p, E_p)$ is the positive (colored in purple) 
and $G_g=(V_g, E_g)$ is the negative (colored in gold) checkerboard graph.
Since $G_p^*=G_g$ we have $|E_p|=|E_g|$ and 
$$|V_p|-|E_p|+|V_g|=2=|V_p|-|E_g|+|V_g|.$$ 

Set $r_p$ (resp. $r_g$) to be the number of vertices in $G_p$ (resp. $G_g$) 
of valency at least $3$.

In \ref{Lackenby} it is proved, using results of Ian Agol that
$$\mbox{Vol}(S^3-K) \geq 2 v_0 (\max (r_p, r_g)-2).$$

If $\tilde G_p=(V_p, \tilde E_p)$ and $\tilde G_g=(V_g, \tilde E_g)$ are the
reduced graphs of $G_p$ and $G_g$ than it is easy to see that

\begin{eqnarray*}
r_p&=& |V_p| - (|E_g|- |\tilde E_g|)\\
&=& 2 - |V_g|+|\tilde E_g| \\
&=& |a_{n+1}|+1.\\
\end{eqnarray*}

Similarily, $r_g=|a_{m-1}|+1$ and the lower bound follows.
\end{proof} 

\subsection {Example}

\begin{figure}
\centerline{{\scalebox{0.38}
{\includegraphics{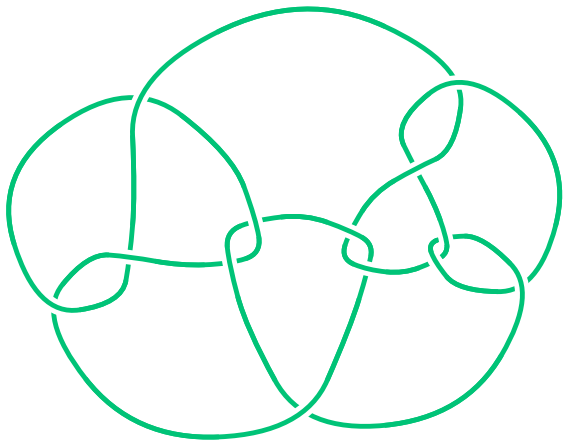}}} \label{SampleKnot}
}
\caption{The alternating knot 13.123 in the Knotscape Census} 
\end{figure}

The checkerboard graph $G$ of the knot in Figure \ref{SampleKnot}
has $|V|=8$ vertices, $|\tilde E|=11, n(2)=2$ and $tri=1$.

It's dual has $|V^*|=7$ vertices, $|\tilde {E^*}|=10, n^*(2)=3$ and $tri^*=2$. 

Therefore, with the above notation for the coefficients of the Jones polynomial,
\begin{eqnarray*}
|a_n|&=&1\\ 
|a_{n+1}|&=&|\tilde E|+1-|V|=4\\
|a_{n+2}|&=&{|a_{n+1}|+1 \choose 2} + n(2) -tri = 10 + 2 -1 = 11\\
|a_m|&=&1\\
|a_{m-1}|&=& |\tilde{E^*}|+1-|V^*|=4\\
|a_{m-2}|&=& {|a_{m-1}|+1 \choose 2} + n^*(2) -tri^* = 10 + 3 -2 = 11
\end{eqnarray*}

The complete Jones polynomial of the knot is, according to Knotscape:
\begin{eqnarray*}
V_{13.121}(t)&=& t^{-12}-4 t^{-11}+11 t^{-10} - 23 t^{-9}+35 t^{-8}-47 t^{-7} +53 t^{-6}\\
&& -52 t^{-5}+47 t^{-4}-34 t^{-3}+22 t^{-2}-11 t^{-1} +4 - t
\end{eqnarray*}
and the hyperbolic volume is $$\mbox{Vol}(S^3-K) \approx 21.1052106828.$$

\newpage
\begin{appendix}
\section{Higher twist numbers of prime alternating knots on $14$ crossings}

In this appendix we give experimental data on the relationsship of the twist number, as computed using the Jones polynomial, and the hyperbolic volume of 
knots. All data are taken from Knotscape, written by Jim Hoste, 
Morwen Thistlethwaite and Jeff Weeks \cite{HTW:Knottabulation}.
We confined ourselves to knots with crossing number 14. 
As before, let $V_K(t)=a_n t^n + a_{n+1} t^{n+1} + \dots + a_m t^m$ 
be the Jones polynomial of an alternating prime knot $K$.

As shown, the twist number is $T(K)=|a_{n+1}|+|a_{m-1}|$.
We call $T_i(K)=|a_{n+i}|+|a_{m-i}|$ the higher ``twist'' numbers.
In particular, $T(L)=T_1(L)$.  

\begin{figure}[h]
\centerline{{\scalebox{0.5}
{\includegraphics{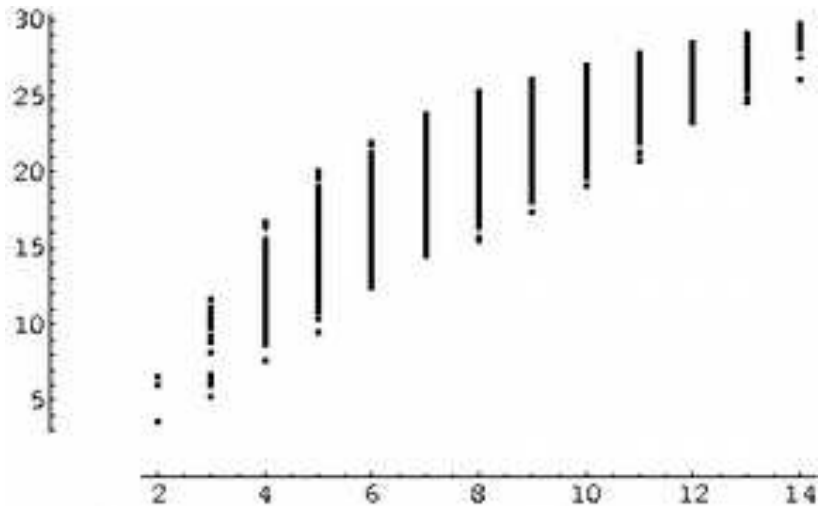}} }}
\caption{The twist number vs. the Volume: 14 crossings, alternating} 
\end{figure}

\begin{figure}[h]
\centerline{{\scalebox{0.5}
{\includegraphics{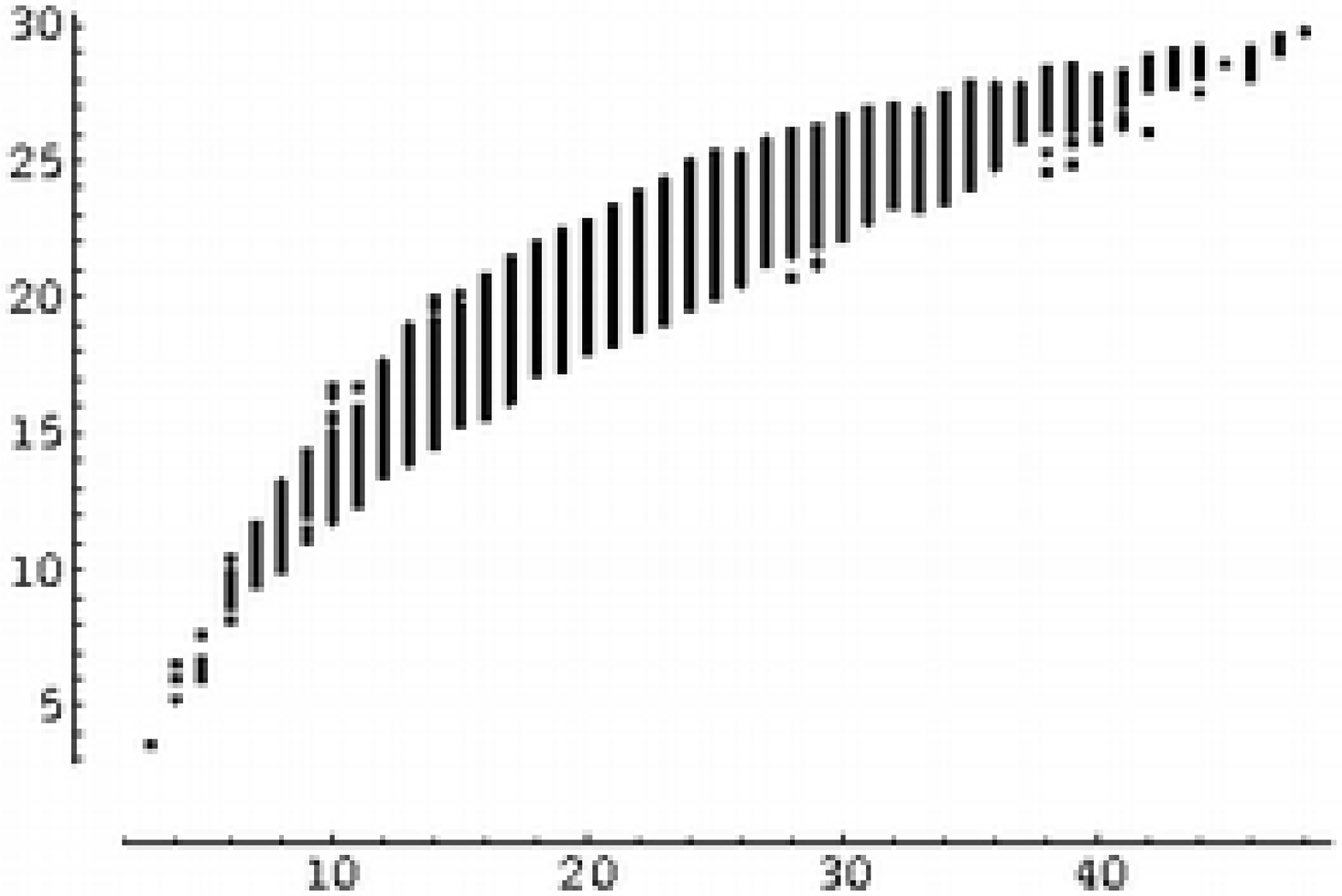}}}
}
\caption{The twist number $T_2$ vs. the Volume: 14 crossings, alternating} 
\end{figure}

\begin{figure}[h]
\centerline{{\scalebox{0.5}
{\includegraphics{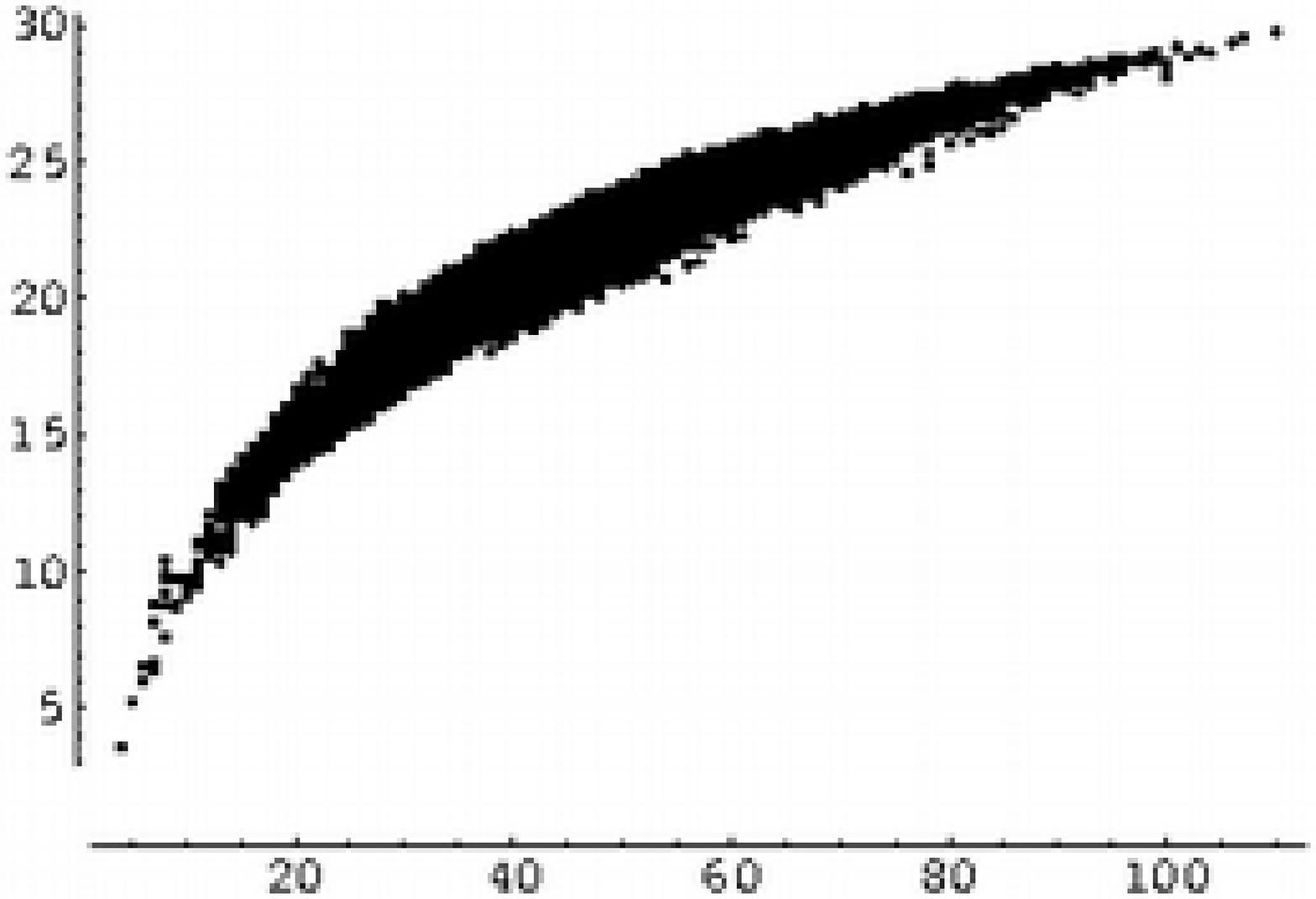}}}
}
\caption{The twist number $T_3$ vs. the Volume: 14 crossings, alternating} 
\end{figure}

\begin{figure}[h]
\centerline{{\scalebox{0.5}
{\includegraphics{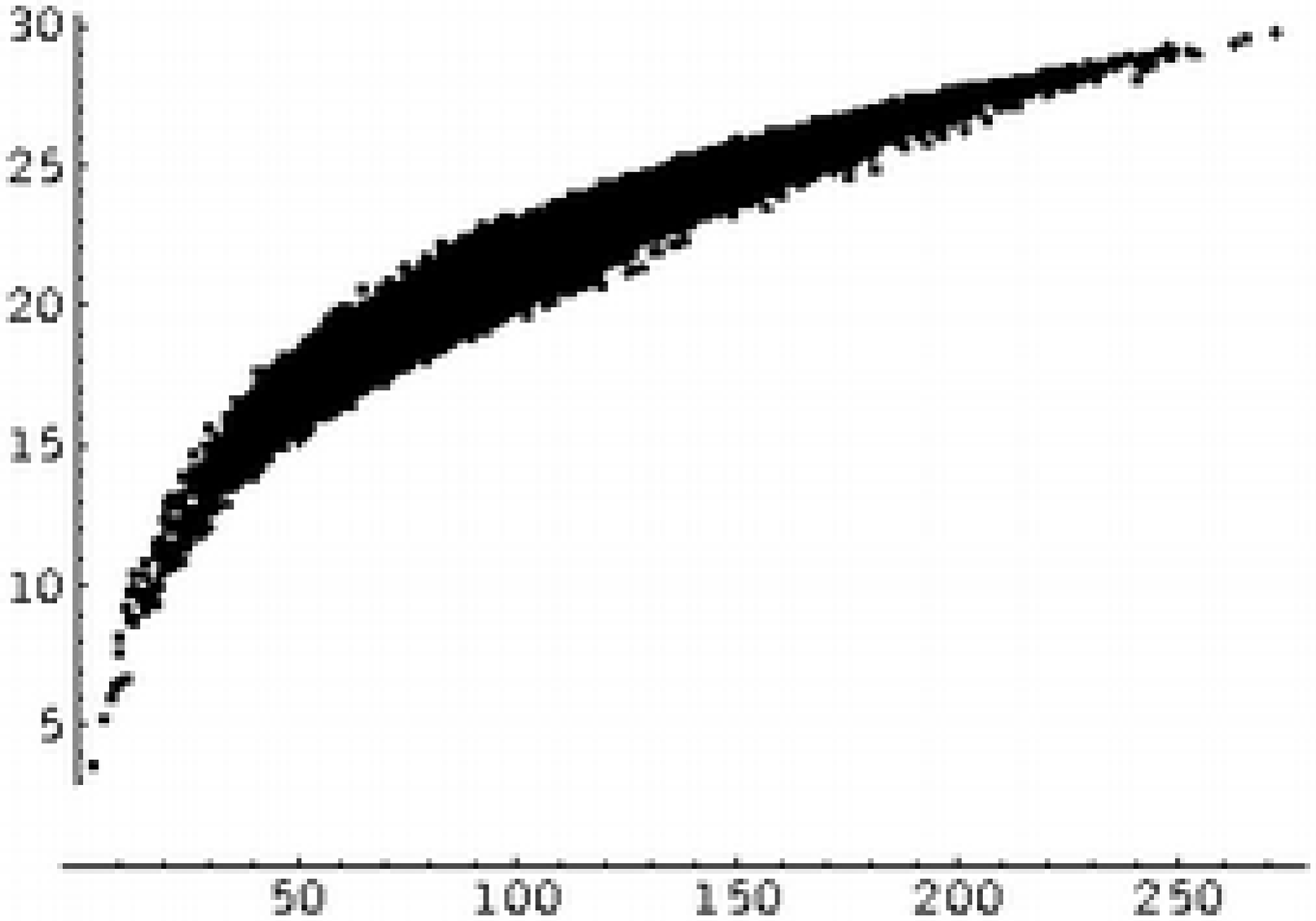}}}
}
\caption{The twist number $T_5$ vs. the Volume: 14 crossings, alternating} 
\end{figure}



\newpage
\section{Higher ``twist numbers'' of prime non-alternating knots on $14$ 
crossings}

For non-alternating knots we keep the notation, 
although there is no direct geometrical justification known:

Again, let $V_L(t)=a_n t^n + a_{n+1} t^{n+1} + \dots + a_m t^m$ 
be the Jones polynomial of a non-alternating knot $L$.

Define the twist number is $T(L)=|a_{n+1}|+|a_{m-1}|$.
As in the alternating case, we call $T_i(L)=|a_{n+i}|+|a_{m-1}|$ 
the higher twist numbers. In particular, $T(L)=T_1(L)$.

If the knot is non-hyperbolic we set its volume = 0.

The pictures give, for non-alternating knots with crossing number 14,
the relation between the twist number and the Volume, resp. the number
$T_2, T_3, T_4$ and the volume.

\begin{figure}[h]
\centerline{{\scalebox{0.5}
{\includegraphics{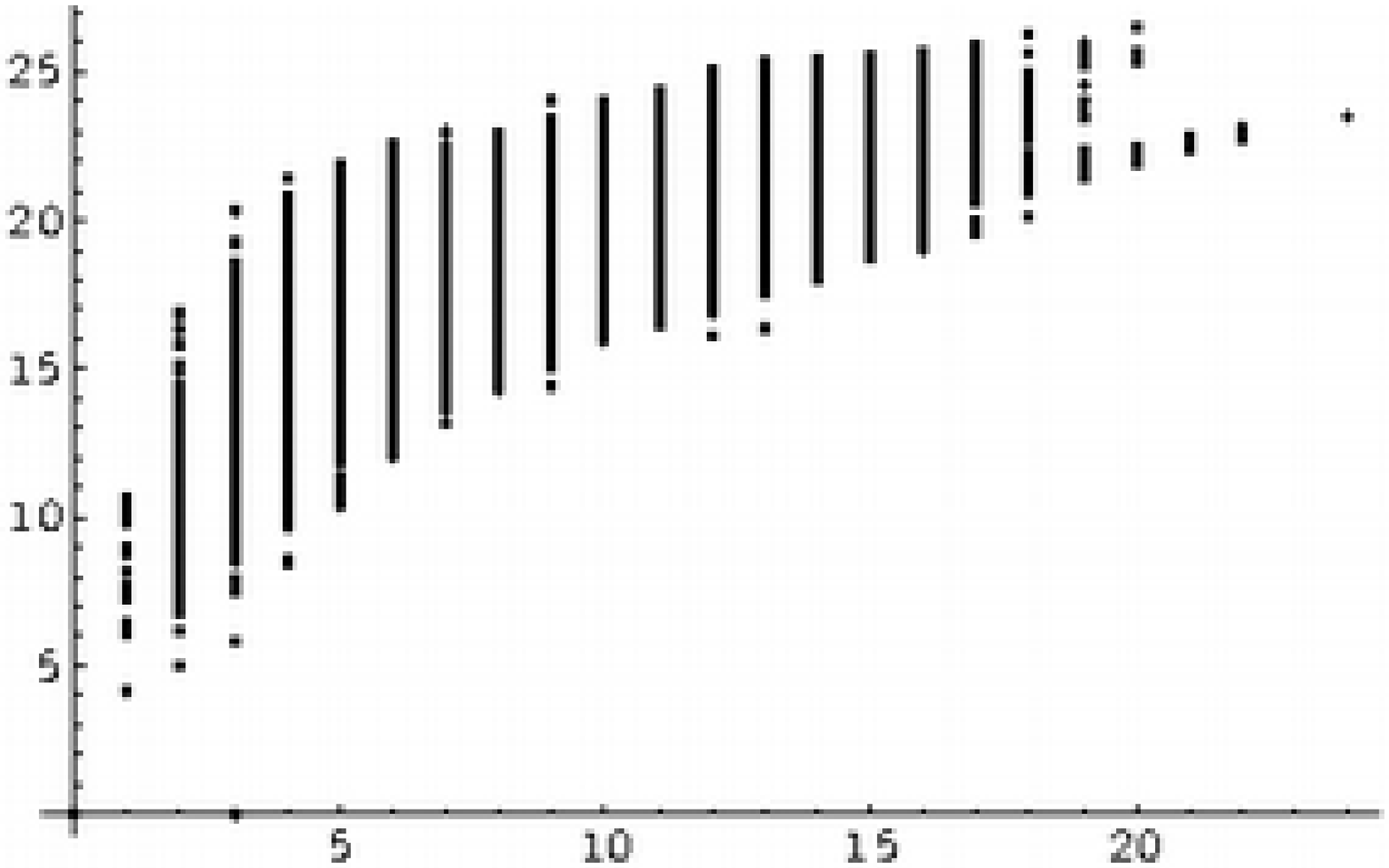}} }}
\caption{The twist number vs. the Volume: 14 crossings, non-alternating} 
\end{figure}

\begin{figure}[h]
\centerline{{\scalebox{0.5}
{\includegraphics{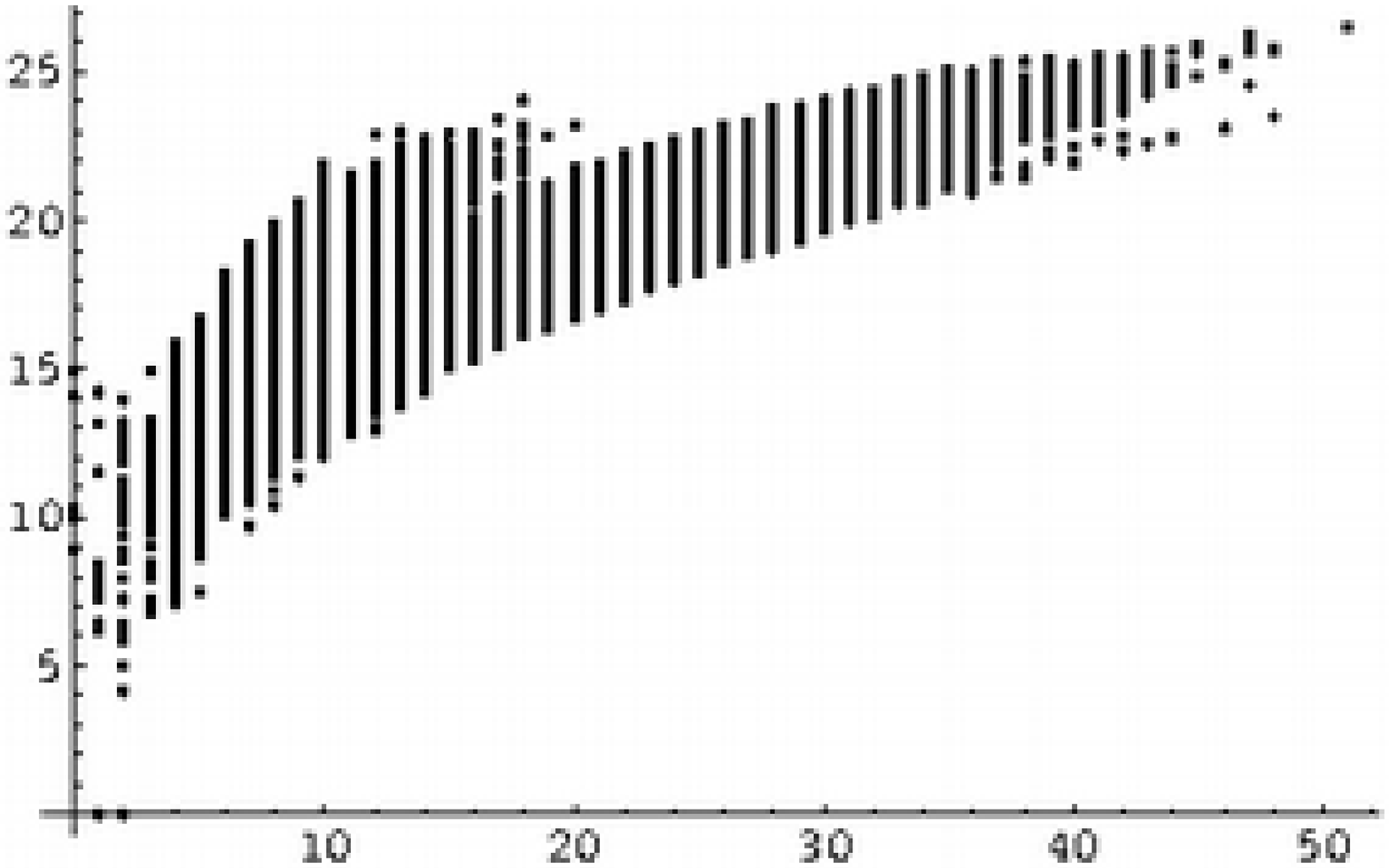}}}
}
\caption{The twist number $T_2$ vs. the Volume: 14 crossings, non-alternating} 
\end{figure}

\newpage

\begin{figure}[h]
\centerline{{\scalebox{0.5}
{\includegraphics{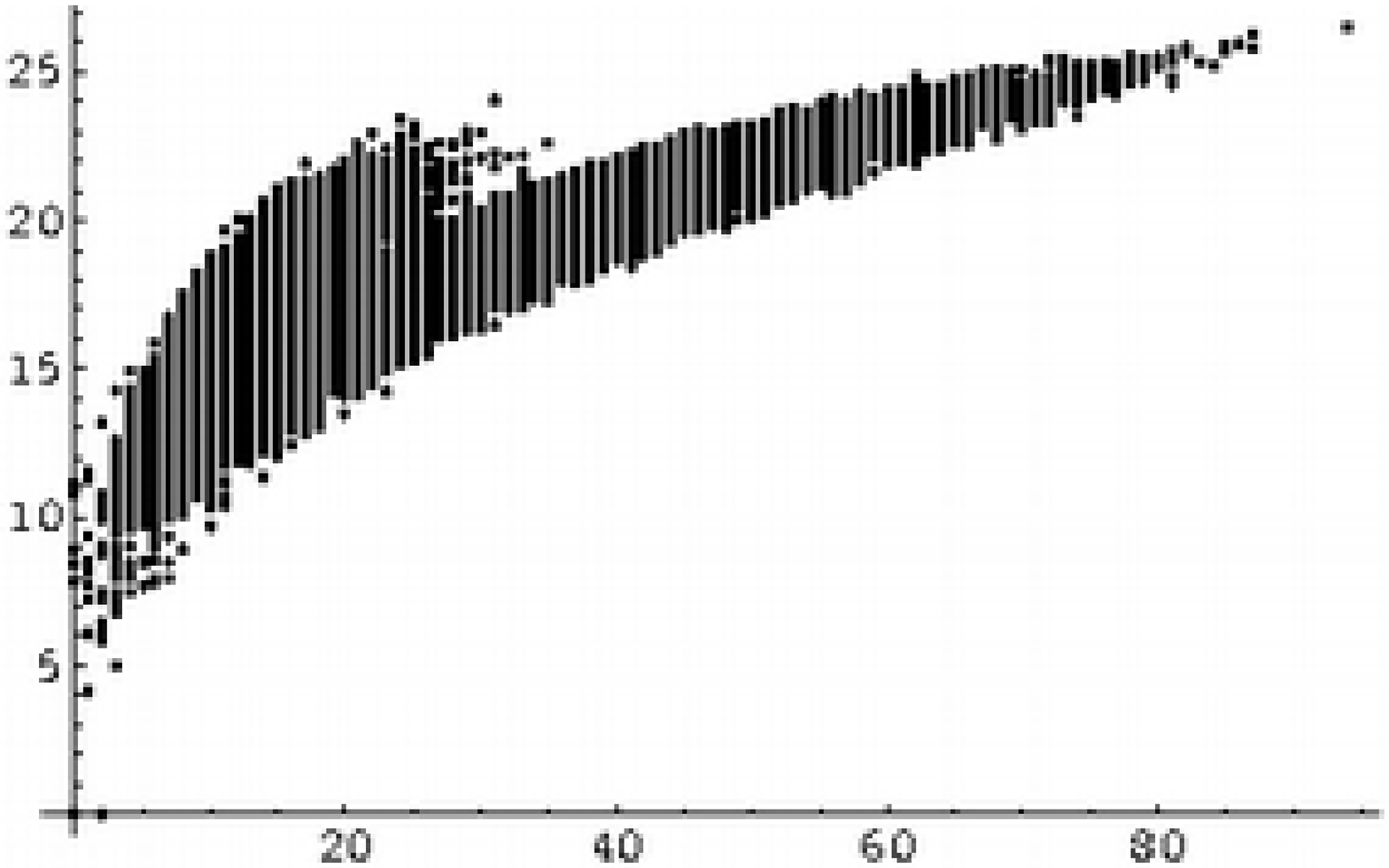}}}
}
\caption{The twist number $T_3$ vs. the Volume: 14 crossings, non-alternating} 
\end{figure}

\begin{figure}[h]
\centerline{{\scalebox{0.5}
{\includegraphics{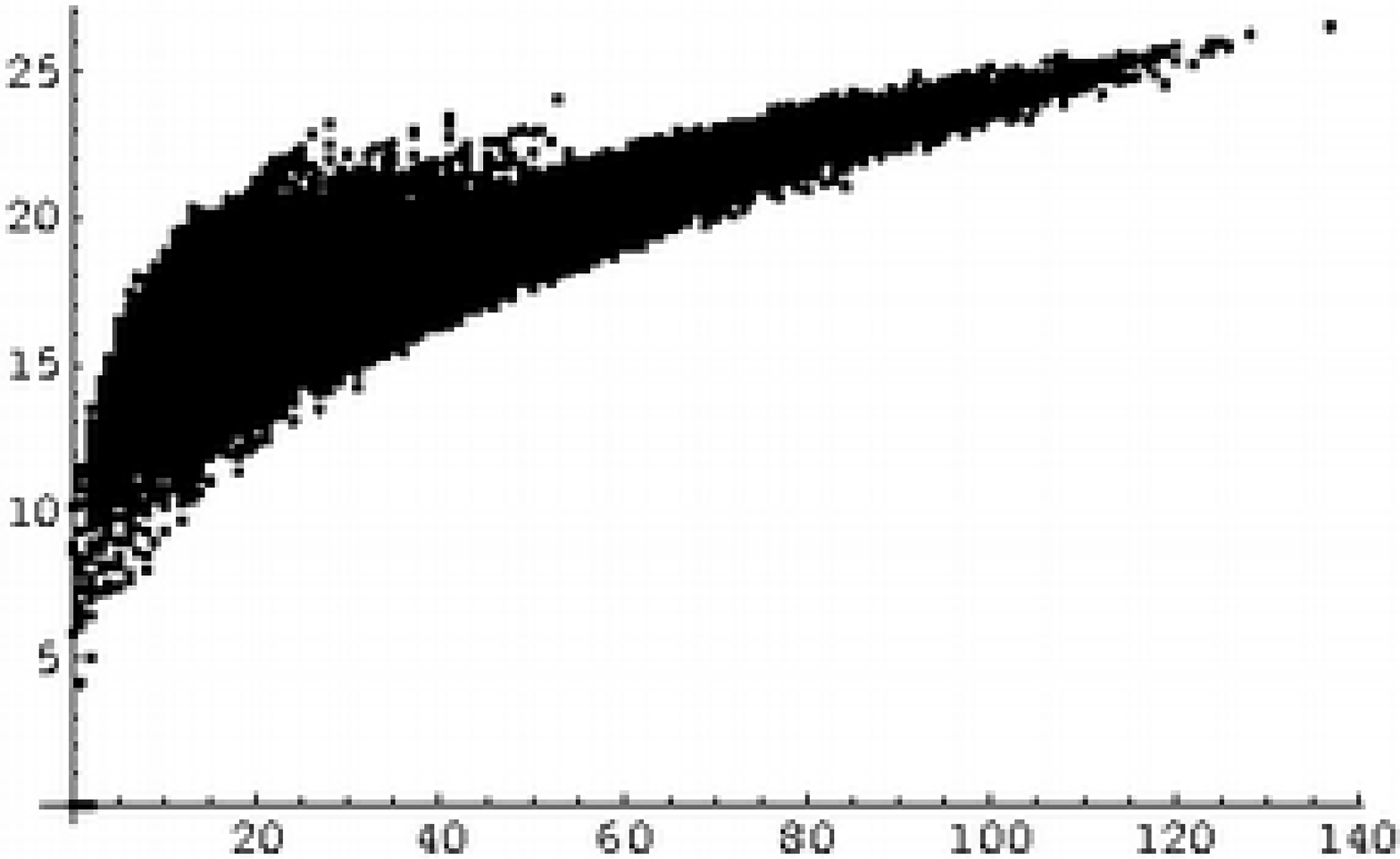}}}
}
\caption{The twist number $T_4$ vs. the Volume: 14 crossings, non-alternating} 
\end{figure}

\end{appendix}


 
 








\bibliography{../linklit}
\bibliographystyle {amsalpha}
\end{document}